\documentclass{article}
\usepackage{shapepar,bm}
\usepackage[dvips,arrow,matrix,ps,color,line]{xy}
\usepackage{theorem,amsfonts,amssymb,amscd}
\usepackage{geometry}
\usepackage{makeidx}
\usepackage{eulervm}

\usepackage{graphicx}
\usepackage{graphics}

\usepackage{epstopdf}

\usepackage{amssymb,graphics,hyperref}

\hypersetup{
    colorlinks = true,
    linkcolor = blue,
    anchorcolor = red,
   citecolor = blue,
    filecolor = red,
    pagecolor = red,
    urlcolor = blue}

\newcommand{\mathsym}[1]{{}}

\usepackage[usenames]{color}
\definecolor{MyLightMagenta}{cmyk}{0.1,0.8,0,0.1}
\definecolor{MyDarkBlue}{rgb}{0.1,0,0.3}

\makeindex

\def\wb{[\bfb]}

\def\NN{\mathbb N}

\def\ovsig{\overline{\sigma}}

\def\bfb{{\mathbf b}}
\def\wi{{w^{-1}}}

\def\gln{{gl_n(\ZZ)}}
\def\glin{{gl_\infty(\ZZ)}}

\def\bfu{{\mathbf u}}
\def\bfv{{\mathbf v}}
\def\bfw{{\mathbf w}}
\def\Qcal{\mathcal Q}

\def\ZZ{\mathbb Z}

\def\CC{\mathbb C}

\def\Ecal{{\mathcal E}}

\def\QQ{\mathbb Q}

\def\cocoa{{\hbox{\rm C\kern-.13em o\kern-.07em C\kern-.13em o\kern-.15em A}}}

\def\Dcal{\mathcal D}

\def\Acal{{\mathcal A}}

\def\hom{{\mathrm{Hom}}}

\def\bfu{{\bf u}}

\def\ovDcal{\overline{\Dcal}}

\def\End{\mathrm{End}}

\def\blamb{{\bm \lambda}}
\def\bmu{{\bm \mu}}

\def\Pcal{{\mathcal P}}

\def\w2M{\bigwedge^2M}

\def\wM{\bigwedge M}

\def\w{\wedge }
\def\bw{\bigwedge }

\protect
\protect
\protect\def\wMn{{\bigwedge M_n}}
\protect

\protect
\protect

\def\sra{\rightarrow}
\def\lra{\longrightarrow}

\def\proof{\noindent{\bf Proof.}\,\,}
\def\qed{{\hfill\vrule height4pt width4pt depth0pt}\medskip}
\def\be{\begin{equation}}
\def\ee{\end{equation}}
\def\bclm{\begin{claim}}
\def\eclm{\end{claim}}
\def\beqn{\begin{eqnarray}}
\def\eeqn{\end{eqnarray}}
\def\beqn*{\begin{eqnarray*}}
\def\eeqn*{\end{eqnarray*}}

\theoremstyle{change}
\theorembodyfont{\rmfamily}

\newtheorem{claim}{}[section]

\global
\global

\def\no@breaks#1{{\def\\{ \ignorespaces}#1}}    

  \makeatletter
\def\cleardoublepage{\clearpage\if@twoside \ifodd\c@page\else
\hbox{} \thispagestyle{empty}
\newpage
\if@twocolumn\hbox{}\newpage\fi\fi\fi} \makeatother

\usepackage{eso-pic,graphicx}
\makeatletter
\newcommand\BackgroundPicture[2]{%
  \setlength{\unitlength}{1pt}%
  default \put(0,\strip@pt\paperheight){%
  \parbox[t][\paperheight]{\paperwidth}{%
    \vfill
     \centering \includegraphics[angle=#2, width=15cm, height=15cm,  bb=0 0 150 150]{#1}
    \vfill
}}} %
\makeatother
\title{The Cohomology of the Grassmannian is a $gl_n$-module}

\author{Letterio Gatto \& Parham Salehyan  \thanks{Work sponsored by FAPESP,  Processo n. 2016/03161-3 and, partially, by Politecnico di Torino, Finanzia-\newline mento 
Diffuso della Ricerca; INDAM-GNSAGA e PRIN "Geometria delle Variet\`a Algebriche''.
\newline ${}$ \,\,\,\,\,\, 2010 MSC: 14M15, 15A75, 17B69.  
\newline ${}$\,\,\,\,\,\,\, {\em Keywords and phrases:} Hasse-Schmidt Derivations on Exterior Algebras, Schubert Derivations;  vertex operators, cohomology of the Grassmannianm bosonic vertex representation of Date-Jimbo-Kashiwara-Miwa.}}
\date{}

\begin{document}
\maketitle
\begin{abstract} 
\noindent 
The integral singular cohomology ring of the  Grassmann variety parametrizing $r$-dimensional subspaces in the $n$-dimensional complex vector space  is naturally an irreducible representation of the Lie algebra of all the $n\times n$ matrices with integral entries.
 Using the notion of Schubert derivation, a distinguished Hasse-Schmidt derivation on an exterior algebra, we describe explicitly such a representation, indicating its relationship with the celebrated bosonic vertex representation of the Lie algebra of infinite matrices due to Date, Jimbo, Kashiwara and Miwa.

\end{abstract}

\section{Introduction}
\claim{} \label{sec:sec01}

%

For  $r,n\in\NN\cup\{\infty\}$ such that $0\leq r\leq n$,  let  $B_{r,n}$ be the singular cohomology ring  of the complex Grassmann variety $G(r,n)$, in the following extended sense. If $r,n$ are both finite, then $B_{r,n}=H^*(G(r,n),\ZZ)$ is the  singular cohomology ring of the usual finite--dimensional Grassmann variety parametrizing $r$-dimensional subspaces of $\CC^n$. If $r<\infty$ and $n=\infty$, then $B_r:=B_{r,\infty}=H^*(G(r,\infty))$, where $G(r,\infty)$ is the ind-variety corresponding to the  chain of inclusions 
$$
\cdots \hookrightarrow G(r, n-1)\hookrightarrow G(r, n)\hookrightarrow G(r,n+1)\hookrightarrow\cdots
$$ In this case $B_r$ is simply a polynomial ring in $r$ indeterminates $\ZZ[e_1,\ldots, e_r]$ (see e.g. \cite{bottu}), which is graded by giving degree $i$ to each indeterminate $e_i$; if both $r,n=\infty$, instead,  $Gr(\infty):=G(r,\infty)$ is the  Universal Grassmann Manifold (UGM) introduced by Sato (see the survey \cite{satoUGM}), and  which is the same as the ind-Grassmannian constructed in \cite{DimiPenkov, IgnaPenkov}. In this case the ring $B:=B_{\infty}$ is the projective limit of $B_r$ in the category of graded modules and, concretely, is a $\ZZ$-polynomial ring $\ZZ[e_1,e_2,\ldots]$ in the infinitely many indeterminates $(e_1,e_2,\ldots)$. 
Let now $\bw  M_n=\bigoplus_{r\geq 0}\bw^rM_n$ be the exterior algebra of the free abelian group  $M_n:=\bigoplus_{0\leq i<n}\ZZ b_i$, with basis $\bfb:=(b_i)_{0\leq i<n}$ and consider the Lie algebra
\be
\begin{tabular}{c}
$\gln:=\{A\in End_\ZZ(M_n)\,|\, Ab_j=0$ for all but finitely many $j\in\NN\}$,
\end{tabular}\label{eq:glnz}
\ee
with respect to the usual commutator.
If $r<\infty$, 
let $\Pcal_{r,n}$ be the set  of all the partitions whose Young diagram is contained in a 
$r\times (n-r)$-rectangle. Then there is a natural $\ZZ$-module isomorphism 
$$
\phi_{r,n}:B_{r,n}\sra 
\bw^rM_n,
$$ 
which maps the basis of Schur polynomials of $B_{r,n}$, parametrized by $\Pcal_{r,n}$ to a natural basis $\wb^r_\blamb$ of $\bw^rM_n$, whose 
elements are labeled by $\Pcal_{r,n}$ as well (see e.g. \cite[Formula (25)]
{pluckercone}). It turns out that $\bw M_n$ is a natural 
representation of $\gln$ and that the exterior powers $\bw^rM_n$ are precisely its 
irreducible sub-representations. The case $r=n=\infty$ needs a few adjustment, 
because the ring $B$ is not isomorphic to any finite exterior power. In this case,   
the latter 
must be replaced by the Fermionic Fock space (FFS),  a suitable irreducible 
representation of a canonical Clifford algebra supported on the direct sum of $
{\frak M}:=\bigoplus_{j\in\ZZ}\ZZ b_j$ with its 
restricted dual. The FFS is naturally a $\glin$-module as 
well. So, in general,
\begin{quotation}
\centerline{ \begin{tabular}{|c|}\hline\cr$B_{r,n}$ is a module over the Lie algebra $\gln$, for all 
$r\leq n\in\NN\cup\{\infty\}$.\cr\cr\hline\end{tabular} }
\end{quotation}
\vspace{-4pt}
  If  $r=\infty$,  its structure is described  by the  
bosonic vertex representation due to Date, Jimbo, Kashiwara 
and Miwa (DJKM) \cite{DJKM01,jimbomiwa}. It amounts to determine the shape of the 
generating function $\Ecal(z,w)=\sum_{i,j\in\ZZ}\Ecal_{i,j}z^iw^{-j}$ of all the 
elementary endomorphisms $\Ecal_{ij}\in \End_\ZZ(\frak M)$, defined by $\Ecal_{ij}
b_k=b_i \delta_{jk}$, acting on $B$. In \cite{GatSal5} we extended the notion of Schubert 
derivation, a priori only defined on an exterior algebra,  to the FFS and, as a 
byproduct, we offered an 
alternative  deduction of the DJKM generating function $\Ecal(z,w)$.
The present paper, instead,  is concerned with the description of $\Ecal(z,w)_n:=\sum_{0\leq i,j<n}\Ecal_{ij}\cdot z^iw^{-j}$ as 
acting on $B_{r,n}$ for $r<\infty$. The formula obtained in this case are new and 
their deduction is cheap.

\claim{} To achieve our goal, we first consider the easiest case, namely 
$n=\infty$, and we 
determine two  equivalent, although looking different, expressions for the action of $\Ecal(z,w):=
\sum_{i,j\geq 0}\Ecal_{ij}z^iw^{-j}$ on $B_r$. They have both interesting 
complementarty features. The former
is useful for explicit computations, and can be implemented as well in the case when 
$r=\infty$ (Cf. \cite[Section 7]{GatSal5}), while the latter visibly shows its close 
relationship with the DJKM representation. In fact it is expressed via suitable 
approximations of the bosonic vertex operators whose expression, in terms of the Schubert derivations introduced in Secton~
\ref{sec:sec3},  is exactly the same  occurring in the DJKM case, as shown in \cite{GatSal5}.

Finally, the $\gln$ structure of $B_{r,n}$, for finite $n$,  will be obtained from the $\glin$--structure of $B_r$, by projecting it through the canonical epimorphism $B_r\sra B_{r,n}$. We shall provide a few examples to show  how our formulas work to write explicit expressions for the product $\Ecal_n(z,w)$ with elements of $B_{r,n}$, something that can be done automatically on a computer.


\claim{}\label{sec:sec02} To state more precisely the results of this paper, let us 
introduce some further piece of notation. The canonical $\gln$-module structure of $\bw M_n$  
($n\in\NN\cup\{\infty\}$) is defined by mapping all $A\in \gln$ to $\delta(A)\in 
\End_Z(\bw M_n)$ such that:
\be
\left\{\matrix{\delta(A)\bfu&=&A\cdot \bfu, &\forall \bfu\in M_n=\bw^1M_n,\cr\cr
\delta(A)(\bfv\w \bfw)&=&\delta(A) \bfv\w \bfw+\bfv\w \delta(A)\bfw,&\forall 
\bfv,\bfw\in \bw M_n}.\right.\label{eq0:repgln}
\ee
Since all $\bfu\in\wMn$ is a finite linear combination of monomials of given 
degree, the initial condition and the Leibniz rule determine the map $\delta(A)$ 
over all $\bw M_n$.  An easy check shows that the commutator $[\delta(A),
\delta(B)]\in \End_\ZZ(\bw M_n)$ is equal to $\delta([A,B])$, where $[A,B]$ 
is the commutator in $\gln$.  The composition of  $\delta$ with the 
restriction map to one 
degree of the exterior algebra,  $A\mapsto\delta(A)_{|\bw^rM_n}$,  turns $
\bw^rM_n$ into a representation of $\gln$. This is easily seen to be irreducible, 
because any basis element of $\bw^rM_n$ can be transported to any other via a 
suitable element of $\gln$. Let us set $M:=M_\infty$. In the ring $B_r:=\ZZ[e_1,\ldots, e_r]$ consider the generic polynomial of degree $r$:
$$ 
E_r(z):=1-e_1z+\cdots+(-1)^re_rz^r\in B_r[z]
$$
and the sequence  $H_r:=(h_j)_{j\in\ZZ}$  implicitly defined by
$$
\sum_{i\geq 0}h_iz^i:={1\over E_r(z)}\in B_r[[z]].
$$ 
It is well known that $B_r=\bigoplus_{\blamb\in\Pcal_r}\ZZ\Delta_\blamb(H_r)$, 
where $\Pcal_r$ denotes the set of all partitions of length at most $r$ and $
\Delta_\blamb(H_r)$ is the Schur determinant $\det(h_{\lambda_j-j+i})_{1\leq i,j< n}$. 
We have a $\ZZ$-module isomorphism $\phi_r:B_r\mapsto \bw^rM$ given by $\Delta_
\blamb(H_r)\mapsto \wb^r_\blamb$, where
$$
\wb^r_\blamb:=b_{r-1+\lambda_1}\w\cdots\w b_{\lambda_r}\in \bw^rM.
$$
Consider now the following two sequences of elements of the polynomial ring $B_r[z^{-1}]$
$$
\ovsig_-(z)H_r:=(\ovsig_-(z)h_j)_{j\in\ZZ}\qquad \mathrm{and}\qquad \sigma_-(z)H_r:=(\sigma_-(z)h_j)_{j\in\ZZ},
$$
where
$$
\ovsig_-(z)h_j=h_j-{h_{j-1}\over z}\qquad \mathrm{and}\qquad \sigma_-(z)h_j=\sum_{i\geq 0}{h_{j-i}\over z^i}.
$$
Extending the action of $\sigma_-(z)$ and $\ovsig_-(z)$ to all basis elemet of $B_r$ through the rule
\be
\sigma_-(z)\Delta_\blamb(H_r)=\Delta_\blamb(\sigma_-(z)H_r)\qquad \mathrm{and}\qquad \ovsig_-(z)\Delta_\blamb(H_r)=\Delta_\blamb(\ovsig_-(z)H_r),
\ee
provides two well defined $\ZZ$-linear maps $\ovsig_-(z),\sigma_-(z):B_r\sra B_r[z^{-1}]$. 
Let  now 
$$
\delta(z,w):=\displaystyle{\sum_{i,j\geq 0}}\delta(\Ecal_{ij})z^iw^{-j}
$$
 and define $\Ecal(z,w):=\displaystyle{\sum_{i,j\geq 0}}b_i\otimes\beta_j\cdot z^iw^{-j}:B_r\sra B_r[[z,w^{-1}]$ via  the equality
 $$
 \phi_{r}(\Ecal(z,w)P(e_1,\ldots,e_r))=\delta(z,w)\phi_{r}(P(e_1,\ldots,e_r)),$$
 where $P(e_1,\ldots,e_r)$ is an arbitrary $\ZZ$-polynomial in $e_1,\ldots,e_r$.
 
 \smallskip
 We are now in position to state the first main result of this paper.

\smallskip
\noindent
{\bf Theorem \ref{thm:thm43}} {\em The $\Ecal(z,w)$-image in $B_r[z,w^{-1}]$ of a basis element $\Delta_\blamb(H_r)\in B_r$, is:
\be
\Ecal(z,w)\Delta_\blamb(H_{r})={z^{r-1}\over w^{r-1}}\cdot  {1\over E_r(z)}\Delta_\blamb(w^{-\blamb},\ovsig_-(z)H_r), \label{eq0:1stmain}
\ee

\medskip
where
\be
\Delta_\blamb(w^{-\blamb},\ovsig_-(z)H_r):=\left|\matrix{w^{-\lambda_1}&w^{1-\lambda_2}&\cdots&w^{r-1-\lambda_r}\cr\cr
h_{\lambda_1+1}-\displaystyle{h_{\lambda_1}\over z}&h_{\lambda_2}-\displaystyle{h_{\lambda_2-1}\over z}&\cdots&h_{\lambda_r+r-2}-\displaystyle{h_{\lambda_r+r-3}\over z}\cr \vdots&\vdots&\ddots&\vdots\cr\cr
h_{\lambda_1+r-1}-\displaystyle{h_{\lambda_{1}+r-2}\over z}&h_{\lambda_2+r-2}-\displaystyle{h_{\lambda_2+r-3}\over z}&\cdots&h_{\lambda_r}-\displaystyle{h_{\lambda_r-1}\over z},
}\right|.\label{eq:fddjkm}
\ee
}

In other words, the product $\Ecal_{ij}\cdot\Delta_\blamb(H_r)$ is the coefficient of $z^iw^{-j}$ in the expansion of the right--hand side of (\ref{eq0:1stmain}).
A second equivalent expression for the action of $\Ecal(z,w)$ is provided by:

\smallskip
\noindent
{\bf Theorem \ref{thm5:2ndthm}} {\em The following formula holds
\be 
\Ecal(z,w)=\left(1-{z\over w}\right)^{-1}\left({z^r\over w^r}\,\Gamma_r(z,w)-1\right),\label{eq0:scndmnres}
\ee
}
where
$$
\Gamma_r(z,w):={E_r(w)\over E_r(z)}\ovsig_-(z)\sigma_-(w).
$$
\smallskip
Equation (\ref{eq0:scndmnres}) recalls the shape of the bosonic representation of the Lie algebra $\Acal_\infty$ of all the matrices with finitely many non-zero diagonals: see \cite{DJKM01,jimbomiwa,KacRaRoz} and, from  now  on, also \cite[Section 9]{GatSal5}.

\claim{} Suppose now that $n<\infty$ and let: 
$$
B_{r,n}:={B_r\over {(h_{n-r+1},\ldots, h_n})}.
$$
Denote by $\pi_{r,n}:B_r\sra B_{r,n}$ the canonical epimorphism. Then
$$
B_{r,n}=\bigoplus_{\blamb\in\Pcal_r}\ZZ\cdot \pi_{r,n}\Delta_\blamb(H_r)= \bigoplus_{\blamb\in\Pcal_r}\ZZ\cdot \Delta_\blamb(H_{r,n}),
$$
where $H_{r,n}=\pi_{r,n}H_r=(\pi_{r,n}(h_j))_{j\in\ZZ}=(1=h_0,h_1,h_2,\ldots, h_{n-r})$. 
 As remakerd,  $B_{r,n}$ is a $gl_n(\ZZ)$-module. We have
 
 \smallskip
 \noindent
{\bf Theorem \ref{thm6:mnthm}.} {\em The following equality holds in $B_{r,n}[z,w^{-1}]$:
\be
 \Ecal(z,w)_n\Delta_{\blamb}(H_{r,n})=\pi_{r,n}(\Ecal(z,w)\Delta_\blamb(H_r)),\label{eq00:proje}
\ee
where $\Ecal(z,w)_n=\sum_{0\leq i,j< n}\Ecal_{i,j}z^iw^{-j}$.
 }

\smallskip
Equality~(\ref{eq00:proje}) means that to describe the $\gln$-action on an element of $B_{r,n}$ is sufficient to set to zero all the $h_j$ with $j>n-r$ that may possibly occur in the expression obtained for $n=\infty$ using (\ref{eq0:1stmain}).
For example, by applying the described recipe, it is easy to see that
\begin{eqnarray*}
\Ecal_{4}(z,w)\Delta_{(2,2)}(H_{2,4})&=&{1\over w^2}\left(-h_2-h_1h_2z+h_2^2z^2\right)+{1\over w^3}\left( -h_1-(h_1^2-h_2)z+h_2^2z^3\right)
\cr\cr
&=&{1\over w^2}[(e_2-e_1^2)+(e_1e_2-e_1^3)z+ (e_1^4-2e_1^2e_2+e_2^2)z^2]\cr
&-&[e_1+e_2z-(e_1^4-2e_1^2e_2+e_2^2)z^3]\in B_{4}[[z,w^{-1}],
\end{eqnarray*}
a result that can be obviously checked by hand, due to the low values of $r$ and $n$ in this specific example.

\medskip

\noindent{\bf Acknowledgments.} This work started during the stay of the first 
author at the Department of Mathematics of UNESP, S\~ao Jos\'e do Rio Preto, under 
the auspices of FAPESP, Processo n. 2016/03161--3. It continued under the partial 
support of INDAM-GNSAGA and the PRIN ``Geometria delle Variet\`a Algebriche''. A 
short visit of the first author to the second one was supported by the program 
``Finanziamento Diffuso della Ricerca'' of Politecnico di Torino. All these 
institutions are warmly acknowledged. For discussions and criticisms we want to 
thank primarily Joachim Kock, as well as, in alphabetical order,  Simon G. Chiossi, 
Abramo Hefez, Marcos Jardim,  Daniel Levcovitz, Simone Marchesi, Igor Mencattini, Piotr Pragacz, Andrea T. Ricolfi, Inna Scherbak and Aron Simis.

\section{\bf Preliminaries and Notation} \label{sec:sec1}
\claim{} A partition is a monotonic non increasing sequence of non-negative integers $\lambda_1\geq \lambda_2\geq\ldots$ all zero but finitely many, said to be {\em parts}. The length $\ell(\blamb)$ of a partition $\blamb$ is the number of non zero parts.  We denote by $\Pcal$ the set of all partitions, by $\Pcal_r:=\{\blamb\in\Pcal_r\,|\,\ell(\blamb)\leq r \}$ and by $\Pcal_{r,n}$ the set of all partitions of length at most $r$ whose Young diagram is contained in a $r\times (n-r)$ rectangle. The partitions form an additive semigroup: if $\blamb,\bmu\in\Pcal$, then $\blamb+\bmu\in\Pcal$.
If $\blamb:=(\lambda_1,\lambda_2,\ldots)$, we denote by $\blamb^{(i)}$ the partition obtained by removing the $i$-th part:
$$
\blamb^{(i)}:=(\lambda_1\geq\lambda_{i-1}\geq\widehat{\lambda_i}\geq\lambda_{i+1}\geq\ldots),
$$
where \,\, $\widehat{}$\,\, means removed. 
By $(1^j)$ we mean the partition with $j$ parts equal to $1$.

\claim{}  In the following $M$ will denote the free abelian group $\bigoplus_{i\geq 0} \ZZ\cdot b_i$ with basis $\bfb:=(b_i)_{i\geq 0}$. 
  For all partitions $\blamb\in\Pcal_r$, let
\be
\wb^r_\blamb:=b_{r-1+\lambda_1}\w\cdots\w b_{\lambda_r}.
\ee 
Clearly $\bw^rM:=\bigoplus_{\blamb\in\Pcal_r} \ZZ\wb^r_\blamb$.
The {\em restricted dual} of $M$ is $M^*:=\bigoplus_{i\geq 0}\ZZ\beta_i$, where $\beta_i\in \hom_\ZZ(M, \ZZ)$ is such that $\beta_i(b_j)=\delta_{ij}$.  
There is a natural well known identification between $(\bw M_n)^*$ and $\bw M_n^*$, see e.g. \cite[Section 2.6]{pluckercone}.

\claim{\bf Hasse-Schmidt derivations on $\wM$.}  Let $z$ denote an arbitrary formal variable.
A {\em Hasse-Schmidt derivation} (HS) \cite{GatSal4} on $\bw M$ is a homomorphisms of abelian groups $\Dcal(z):\bw M\sra \bw M[[z]]$ such that
\be
\Dcal(z)(\bfu\w \bfv)=\Dcal(z)\bfu\w \Dcal(z)\bfv,\label{eq:defHs}
\ee
for all $\bfu,\bfv\in \bw M$. Writing $\Dcal(z)$ as $\sum_{j\geq 0}D_jz^j$, equation~(\ref{eq:defHs}) is equivalent to
$$
D_j(\bfu\w \bfv)=\sum_{i=0}^jD_i\bfu\w D_{j-i}\bfv.
$$

\claim{}\label{sec:dbar} If $\Dcal(z)$ is a HS--derivation on $\wM$  and  $D_0$ is invertible,  there exists $\ovDcal(z):=\sum_{i\geq 0}(-1)^i\ovDcal_iz^i\in \End_\ZZ(\wM)[[z]]$ such that $\ovDcal(z)\Dcal(z)=\Dcal(z)\ovDcal(z)=1$. The map  $\ovDcal(z)$ is a HS-derivations, said to be the {\em inverse} of $\Dcal(z)$.
Thus the two integration by parts formulas hold:
\be
\Dcal(z)\bfu\w \bfv=\Dcal(z)(\bfu\w \ovDcal(z)\bfv) \qquad\mathrm{and}\qquad \bfu\w \ovDcal(z)\bfv=\ovDcal(z)(\Dcal(z)\bfu\w \bfv).\label{eq:ipart}
\ee
As remarked in \cite{GatSche03}, the second of (\ref{eq:ipart}) is the generalization (holding also for free $A$-module of infinite rank) of the Cayley-Hamilton theorem. 

\claim{\bf The transposed HS--derivation.} For all $\eta\in \bw M^*$, let $\Dcal^T(z)\eta$ be the   unique element of $\bw M^*$ such that
$$
\Dcal^T(z)(\eta)(\bfu)=\eta (\Dcal(z)(\bfu)),
$$
for all $\bfu\in\wM$. By \cite[Proposition 3.8]{pluckercone},  $\Dcal^T(z)$ is a HS derivation on $\wM^*$ said to be {\em the transposed} of $\Dcal(z)$. Integration by parts (\ref{eq:ipart}) implies the following equality for transposed HS--derivations.
\bclm{\bf Proposition.} {\em For all $\eta\in M^*$ and each $\bfu\in \bw^rM$
\be
\Dcal^T(z)\eta\lrcorner \bfu=\ovDcal(z)\left(\eta \lrcorner \Dcal(z)\bfu)\right)\label{eq2:iptra}
\ee
}
\eclm
\proof
By definition of contraction of an exterior vector against a linear form, for all $\zeta\in \bw^{r-1}M^*$:
$$
\zeta(\Dcal(z)^T\eta\lrcorner \bfu)=(\Dcal(z)^T\eta\w \zeta)(\bfu).
$$
Now we apply the first of integration by parts (\ref{eq:ipart}):
$$
(\Dcal(z)^T\eta\w \zeta)(\bfu)=\Dcal(z)^T(\eta\w\ovDcal(z)^T\zeta)(\bfu) 
$$
from which, by definition of transposistion,
$$
(\eta\w \ovDcal(z)^T\zeta)\Dcal(z)\bfu=\ovDcal(z)^T\zeta(\eta\lrcorner \Dcal(z)\bfu)=\zeta\left[\ovDcal(z)(\eta\lrcorner \Dcal(z)\bfu)\right]
$$
which proves (\ref{eq2:iptra}).\qed

\section{Schubert derivations on $\wM$}\label{sec:sec2}

It is easy to check that a $HS$-derivation on $\bw M$ is uniquely determined by its restriction to the first degree $M=\bw^1M$ of the exterior algebra (Cf. \cite[Ch. 4]{GatSal4}). Let
$\sigma_+(z):=\sum_{i\geq 0} \sigma_iz^i:\bw M_n\sra \bw M_n[[z]]$ and $\sigma_-(z):=\sum \sigma_{-i}z^{-i}:\bw M\sra \bw M[[z^{-1}]]$ be the unique HS-derivations  such that for all $i\in\ZZ$,  $\sigma_ib_j=b_{i+j}$ if $ i+j\geq 0$ and $0$ otherwise. Let $\ovsig_+(z)$ and $\ovsig_-(z)$ be, respectively, the inverse HS-derivations of $\sigma_+(z)$ and $\sigma_-(z)$ in the algebra $\End_\ZZ(\wMn)[[z^{\pm 1}]]$.

\bclm{\bf Definition.} {\em The HS-derivations $\sigma_{\pm}(z)$ and $\ovsig_{\pm}(z)$ are called {\em Schubert derivations}.
}
\eclm

Let $\Delta_\blamb(\sigma_+)=\det(\sigma_{\lambda_j-j+i})_{1\leq i,j\leq r}$. 
Giambelli's formula for Schubert derivations \cite[Corollary 5.8.2]{GatSal4} or~\cite[Formula (3.2)]{pluckercone} says that $\wb^r_\blamb=\Delta_\blamb(\sigma_+)\wb^r_0$.
 It enables us to equip $\bw^rM$ with a structure of $B_{r}$-module by declaring that
\be
 h_i\wb^r_\blamb:=\sigma_i\wb^r_\blamb.\label{eq:Bmodst}
\ee
Thus  $\bw^r M$  can be thought of as a free $B_r$-module of rank $1$ generated by $\wb^r_0$,  such that $\wb^r_\blamb=\Delta_\blamb(H_{r})\wb^r_0$. In particular 
\be
\sigma_+(z)\wb^r_\blamb={1\over E_r(z)}\wb^r_\blamb.\label{eq3:sez}
\ee

\claim{} Using the $B_r$-module structure of $\bw^rM$ one can define  $\ovsig_-(z)$ and $\sigma_-(z)$ as maps $B_r\sra B_r[z^{-1}]$, by setting
\[
(\sigma_-(z)\Delta_\blamb(H_r))\wb^r_0=\sigma_-(z)\wb^r_\blamb\qquad \mathrm{and}\qquad (\ovsig_-(z)\Delta_\blamb(H_r))\wb^r_0=\ovsig_-(z)\wb^r_\blamb.
\]
A simple application of the definition shows, as in \cite[Proposition 5.3]{pluckercone}, that the equalities 
\be
\sigma_-(z)h_j=\sum_{i=0}^j{h_{j-i}\over z^i}\qquad \mathrm{and}\qquad \ovsig_-(z)h_j=h_j-{h_{j-1}\over z},
\ee
hold  in $B_r$ for all $r\geq 1$.
\bclm{\bf Proposition.}\label{prop:33} {\em Let $\sigma_-(z)H_r$ {\em (}resp. $\ovsig(z)H_r${\em )} stands for the sequence $(\sigma_-(z)h_j)_{j\in\ZZ}$  {\em (}resp. $(\ovsig_-(z)h_j)_{j\in\ZZ}${\em )}. If $\ell(\blamb)\leq r$ then
$$
\sigma_-(z)\Delta_\blamb(H_r)=\Delta_\blamb(\sigma_-(z)H_r)\qquad \mathrm{and}\qquad \ovsig_-(z)\Delta_\blamb(H_r)=\Delta_\blamb(\ovsig_-(z)H_r).
$$
}
\eclm
\proof According to \cite[Theorem 5.7]{pluckercone}, by using a general determinantal formula due to Laksov and Thorup as in \cite[Main Theorem]{LakTho03}.\qed

\bclm{\bf Lemma.} {\em For all $\blamb\in \Pcal_r$ the following equality holds:
\be
b_0\w \ovsig_+(z)\wb^r_\blamb= z^r\ovsig_-(z)\left(\wb^r_{\blamb+(1^r)}\w b_0\right).\label{eq1:lemut}
\ee
}
\eclm
\proof
One argue by induction on $r\geq 1$.  For $r=1$ one has:
\begin{eqnarray*}
b_0\w \ovsig_+(z)b_\lambda&=&b_0\w (b_\lambda-b_{\lambda+1}z)\cr\cr
&=&-b_0\w z(b_{\lambda+1}-b_\lambda z^{-1})\cr\cr
&=&z(b_{\lambda+1}-b_\lambda z^{-1})\w b_0=z\ovsig_-(z)b_{\lambda+1}.
\end{eqnarray*}
Assume (\ref{eq1:lemut}) holds for all $1\leq s\leq r-1$ . Then

\begin{center}
\begin{tabular}{rllr}
$b_0\w \ovsig_+(z)\wb^r_\blamb$&$=$&$b_0\w \ovsig_+(z)\left(b_{r-1+\lambda_1}\w \wb^{r-1}_{\blamb^{(1)}}\right)$&(decomposition of $\wb^r_\blamb$)\cr\cr
&$=$&$b_0\w \ovsig_+(z)b_{r-1+\lambda_1}\w\ovsig_+(z) \wb^{r-1}_{\blamb^{(1)}}$&($\ovsig_+(z)$ is a HS-derivation)\cr\cr
&$=$&$z\ovsig_-(z)b_{r-1+\lambda_1+1}\w (b_0\w \ovsig_+(z) \wb^{r-1}_{\blamb^{(1)}})$&(case $r=1$)\cr\cr
&$=$&$z\ovsig_-(z)b_{r-1+\lambda_1+1}\w z^{r-1}\ovsig_-(z)\wb^{r-1}_{\blamb^{(1)}+(1^{r-1})}$&(inductive hypothesis)\cr\cr
&$=$&$z^r\ovsig_-(z)(\wb^r_{\blamb+(1^r)}\w b_0).$&($\ovsig_-(z)b_0=b_0$ and $\ovsig_-(z)$\cr
&&& is a HS derivation)
\end{tabular}
\end{center}
\qed

\claim{}\label{sec3:35} It is convenient to introduce one more formal variable $w$. Define 
$$
\bfb(z):=\sum_{j\geq 0}b_jz^j\qquad \mathrm{and}\qquad {\bm\beta}(w)=\sum_{j\geq 0}\beta_jw^{-j-1}.
$$
Then $\bfb(z)=\sigma_+(z)b_0$. Moreover ${\bm\beta}(w)=w^{-1}\sigma^T_-(w)\beta_0$. Indeed $(\sigma^T_{-i}\beta_j)(b_k)=\beta_j(\sigma_{-i}b_k)=\beta_j(b_{k-i})=\delta_{i+j,k}=\beta_{j+i}(b_k)$.

Let $\Gamma_{r}(z):B_r\sra B_{r+1}[[z]]$ and $\Gamma^*_{w}(z):B_{r}\sra B_{r-1}[w^{-1}]$ be the operators implicitly defined by:
\begin{eqnarray}
(\Gamma_{r}(z)\Delta_\blamb(H_{r}))\wb^{r+1}_0&=&z^{-r}\bfb(z)\w \wb^r_\blamb\label{eq:avo1}\\\cr 
(\Gamma^*_{r}(w)\Delta_\blamb(H_{r}))\wb^{r-1}_0&=&w^{r}{\bm\beta}(w)\lrcorner \wb^r_\blamb.\label{eq:avo2}
\end{eqnarray}
Clearly $\Gamma_r(z), \Gamma_r^*(z)$ are the finite $r$ case of the bosonic vertex operators as in \cite{KacRaRoz}.
\bclm{\bf Proposition.} {\em
\be
\Gamma_r(z)\Delta_\blamb(H_r)={1\over E_{r+1}(z)}\Delta_\blamb(\ovsig_-(z)H_{r+1})\label{eq:vo01}
\ee
and
\be
\Gamma_r^*(w)\Delta_\blamb(H_r)=\Delta_\blamb(w^{-\blamb},H_{r-1}),\label{eq:vo02}
\ee
where the notation 
\be
\Delta_\blamb(w^{-\blamb},H_{r-1}):=\left|\matrix{w^{-\lambda_1}&w^{-\lambda_2+1}&\cdots& w^{-\lambda_r+r-1}\cr
h_{\lambda_1+1}&h_{\lambda_2}&\cdots&h_{\lambda_r+r-2}\cr
\vdots&\vdots&\ddots&\vdots\cr
h_{\lambda_1+r-1}&h_{\lambda_2+r-2}&\cdots& h_{\lambda_r}
}\right|\label{eq:vo03}
\ee
is used
to keep track of the fact that  all the $h_j$ occurring in the determinant (\ref{eq:vo03}) live in the ring  $B_{r-1}$.
}
\eclm
\proof Let us prove (\ref{eq:vo01}) first. One has

\begin{center}
\begin{tabular}{rllr}
$\bfb(z)\w \wb^r_\blamb$&$=$&$\sigma_+(z)b_0\w \wb^r_\blamb$& (definition of $\sigma_+(z)$)\cr\cr
&$=$&$\sigma_+(z)(b_0\w \ovsig_+(z)\wb^r_\blamb)$&(integration by parts)\cr\cr
&$=$&$\sigma_+(z)(z^r\ovsig_-(z)\wb^r_{\blamb+(1^r)}\w b_0)$&(Formula\ref{eq1:lemut})\cr\cr
&$=$&$z^r\sigma_+(z)\ovsig_-(z)(\wb^{r}_{\blamb+(1^r)}\w b_0)$&($\ovsig_-(z)b_0=b_0$ and \cr
&&&$\ovsig_-(z)$ is a HS derivation)\cr\cr

&$=$&$\displaystyle{z^r\over E_{r+1}(z)}\left(\sigma_-(z)\wb^{r+1}_\blamb\right)$&(by the $B_r$-module structure\cr
&&& (\ref{eq:Bmodst}) of $\bw^rM$)\cr\cr
&=&$	\displaystyle{z^r\over E_{r+1}(z)}\left(\sigma_-(z)\Delta_\blamb(H_{r+1})\right)\wb^{r+1}_0$&(definition of $\ovsig_-(z)$\cr &&&  as a map $B_r\sra B_r[z^{-1}]$)\cr\cr
&$=$&$\displaystyle{z^r\over E_{r+1}(z)}\Delta_\blamb(\sigma_-(z)H_{r+1})\wb^{r+1}_0$.&(Proposition \ref{prop:33}).
\end{tabular}
\end{center}
To prove  (\ref{eq:vo02}), instead, the best is acting by direct computation:
$$
{\bm\beta}(w)\lrcorner \left(b_{r-1+\lambda_1}\w\cdots\w b_{\lambda_r}\right)
$$
$$
=w^{-r-\lambda_1}\wb^{r-1}_{\blamb^{(1)}}-
w^{-r+1-\lambda_2}\wb^{r-1}_{\blamb^{(2)}+(1)}+
\cdots+(-1)^{r-1}w^{-\lambda_r-1}\wb^{r-1}_{\lambda^{(r)}+(1^{r-1})}
$$
$$
=w^{-r}\left(w^{-\lambda_1}\Delta_{\blamb^{(1)}}(H_{r-1})-
w^{1-\lambda_2}\Delta_{\blamb^{(1)}}(H_{r-1})+\cdots +(-1)^{r-1}w^{r-1-\lambda_r}\Delta_{\blamb^{(r)}+(1^{r-1})}\right)\wb^{r-1}_0
$$
$$
=w^{-r}\left|\matrix{w^{-\lambda_1}&w^{-\lambda_2+1}&\cdots& w^{-\lambda_r+r-1}\cr
h_{\lambda_1+1}&h_{\lambda_2}&\cdots&h_{\lambda_r+r-2}\cr
\vdots&\vdots&\ddots&\vdots\cr
h_{\lambda_1+r-1}&h_{\lambda_2+r-2}&\cdots& h_{\lambda_r}
}\right|\wb^{r-1}_0=w^{-r}\Delta_\blamb(w^{-\blamb},H_{r-1})\wb^{r-1}_0,
$$
from which the desired expression of $\Gamma^*_r(w)\Delta_\blamb(H_r)$.\qed

\section{The $\glin$ structure of $B_r$. First description} \label{sec:sec3}
Let $\Ecal_{ij}:=b_i\otimes \beta_j\in \glin$. 
\bclm{\bf Proposition.}
$$
\delta(\Ecal_{ij})\wb^r_\blamb=b_i\w (\beta_j\lrcorner \wb^r_\blamb).
$$
\eclm
\proof It is an easy check, provided one invokes the very definition of the contraction operator as a derivation of degree $-1$ of the exterior algebra.\qed

\claim{} Let $\delta(z,w)=\sum_{i,j\geq 0}\delta(\Ecal_{ij})z^iw^{-j}$. Then $
\delta(z,w)\wb^r_\blamb\in \bw^rM[[z,w^{-1}]
$.
Define  $\Ecal(z,w):B_r\sra B_r[[z,w^{-1}]$ through the equality:
$$
(\Ecal(z,w)\Delta_\blamb(H_r))\wb^r_0:=\delta(z,w)\wb^r_\blamb.
$$

\bclm{\bf Theorem.}\label{thm:thm43}
\be
\Ecal(z,w)\Delta_\blamb(H_r)={z^{r-1}\over w^{r-1}}{1\over E_r(z)}\ovsig_-(z)\Delta_\blamb(w^{-\blamb},H_r)={z^{r-1}\over w^{r-1}}{1\over E_r(z)}\Delta_\blamb(w^{-\blamb},\ovsig_-(z)H_r).\label{eq5:1stmnthm}
\ee
\eclm
\proof
Since $\bw^rM$ is a free $B_r$-module generated by $\wb^r_0$, it suffices to expresses $\delta(z,w)\wb^r_\blamb$ as a $B_r[[z,w^{-1}]$-multiple of $\wb^r_0$. One has:
\begin{eqnarray}
\delta(z,w)\wb^r_0&=&\bfb(z)\w (w{\bm\beta}(w)\lrcorner\wb^r_\blamb)\cr\cr
&=&\sigma_+(z)b_0\w (w\cdot w^{-r}\Gamma^*_r(w)\Delta_\blamb(H_r)\wb^{r-1}_0)\cr\cr
&=&w^{-r+1}\sigma_+(z)(b_0\w \ovsig_+(z)\big(\Gamma^*_r(w)\Delta_\blamb(H_r)\big)\wb^{r-1}_0).\label{eq5:plug}
\end{eqnarray}
Since $(\Gamma^*_r(w)\Delta_\blamb(H_r))\wb^{r-1}_0$ is a  finite linear combination $\sum_{\bmu\in \Pcal_{r-1}}a_\bmu(w)\wb^{r-1}_\mu$:
\begin{eqnarray}
b_0\w \ovsig_+(z)\sum_\bmu a_\bmu(w)\wb^{r-1}_\blamb&=&\sum_\bmu a_\bmu(w)(b_0\w \ovsig_+(z)\wb^{r-1}_\blamb)\cr
&=&\sum_\bmu a_\bmu(w)\left(z^{r-1}\ovsig_-(z)\wb^{r-1}_{\bmu+(1^{r-1})}\w b_0\right)\cr
&=&z^{r-1}\ovsig_-(z)\sum a_\bmu(w)\wb^r_\bmu\cr\cr
&=&z^{r-1}\ovsig_-(z)\Delta_\blamb(w^{-\blamb},H_r)\wb^r_0.\label{eq5:toplug}
\end{eqnarray}
Plugging (\ref{eq5:toplug}) into~(\ref{eq5:plug})
one finally obtains the equality
\be
\delta(z,w)\wb^r_\blamb={z^{r-1}\over w^{r-1}}\sigma_+(z)(\ovsig_-(z)\Delta_\blamb(w^{-\blamb}, H_r)\wb^r_0.\label{eq4:repl}
\ee
Using the $B_r$-module structure of $\bw^rM$ over $B_r$, one may replace $\sigma_+(z)$ by $1/E_r(z)$ in~(\ref{eq4:repl}), getting:
\be
(\Ecal(z,w)\Delta_\blamb(H_r))\wb^r_0\delta(z,w)\wb^r_\blamb={z^{r-1}\over w^{r-1}}{1\over E_r(z)}(\ovsig_-(z)\Delta_\blamb(w^{-\blamb}, H_r)\wb^r_0,\label{eq5:eisi}
\ee
from which, by comparing the coefficients of $\wb^r_0$ on either side of (\ref{eq5:eisi}), and using Proposition \ref{prop:33}, precisely (\ref{eq5:1stmnthm}).
\qed
\bclm{\bf Example.}\label{ex:ex44} Let us compute $\Ecal(z,w)e_2$ in $B_2[[z,w^{-1}]$.
Remind that
$$
e_2=\Delta_{(1,1)}(H_2):=\left|\matrix{h_1&1\cr h_2&h_1}\right|
$$
 and corresponds to the basis element $\wb^2_{(1,1)}:=b_2\w b_1\in\bw^2M$, In particular we expect that $\Ecal_{ij}e_2=0$ for all $j\notin\{1,2\}$.
 By applying the recipe:
\begin{eqnarray}
 \Ecal(z,w)e_2&=&{z\over w}(1+h_1z+h_2z^2+\cdots)\left|\matrix{w^{-1}& 1\cr\cr h_2-\displaystyle{h_1\over z}&h_1-\displaystyle{1\over z}}\right|\cr\cr
 &=&
 {1\over w}{h_1-h_2z\over E_2(z)}+{1\over w^2}{h_1z-1\over E_2(z)}\cr\cr
 &=&\left[{1\over w}(h_1-h_2z)+{1\over w^2}(h_1z-1)\right](1+h_1z+h_2z^2+h_3z^3+\cdots)\label{eq5:23}
\end{eqnarray}
So, for instance
\begin{center}
$\Ecal_{4,2} e_2=$ coefficient of $z^4w^{-2}$ of (\ref{eq5:23})$=h_1h_3-h_4=\Delta_{(3,1)}(H_2)=e_1^2e_2-e_2^2$.
\end{center}
\eclm

\section{The $\glin$ structure of $B_r$. Second description} \label{sec:sec4}
We now compute an equivalent expression of the generating function $\Ecal(z,w)$ which recalls the shape of  the bosonic vertex representation of the Lie algebra $\Acal_\infty$ of the matrices of infinite size with only finitlely many non-zero diagonals.
\claim{} Recall that $\bfb$ and ${\bm\beta}$ satisfy the Clifford algebras relations:
$$
b_i\w (\beta_j\lrcorner)+\beta_j\lrcorner (b_i\w)=\delta_{ij}.
$$
Thus
\be
w(\bfb(z) \w {\bm\beta}(w)\lrcorner \wb^r_\blamb+ {\bm\beta}(w)\lrcorner (\bfb(z)\w \wb^r_\blamb)=\sum_{i\geq 0}{z^i\over w^i}\wb^r_\blamb=i_{w,z}{w\over w-z}\wb^r_\blamb,\label{eq5:ww-z}
\ee
where, following \cite[p.~18]{Kacbeg}, the $i_{w,z}$ means that we are considering the  expansion of $w/(w-z)$  in power series of $z/w$.
We can then compute the $\glin$-action of $B_r$ by first computing $w{\bm\beta}(w)\lrcorner (\bfb(z)\w \wb^r_\blamb)$ and then subtracting the last term of equality (\ref{eq5:ww-z}).

\bclm{\bf Lemma.} {\em For all $i\geq 1$:
\be
\sigma_-(w)b_{n+i}=\sigma_{i}\sigma_-(w)b_n+{1\over w^{n+1}}\sigma_-(w)b_{i-1}.\label{eq5:19}
\ee
}
\eclm
\proof
In fact 
\begin{eqnarray*}
\sigma_-(w)b_{n+i}&=&b_{n+i}+{b_{n+i-1}\over w}+\cdots+{b_i\over w^n}+{1\over w^{n+1}}\sigma_-(w)b_{i-1}\cr \cr
&=&\sigma_i\sigma_-(w)b_n+{1\over w^{n+1}}\sigma_-(w)b_{i-1},
\end{eqnarray*}
as desired.\qed

\bclm{\bf Lemma.}\label{lem:lem56} {\em The following commutation rule holds:
$$
\sigma_-(w)\sigma_+(z)b_0=i_{w,z}{w\over w-z}\sigma_+(z)\ovsig_-(w)b_0.
$$

}
\eclm
\proof
Indeed
$$
\sigma_-(w)\sigma_+(z)b_0=\sigma_-(w)(b_0+b_1z+b_2z^2+b_3z^3+\cdots)
$$
$$
=b_0+ \left({b_0\over w}+b_1\right)z+ \left({b_0\over w^2}+{b_1\over w}+b_2\right)z^2+\left({b_0\over w^3}+{b_1\over w^2}+{b_2\over w}+b_3\right)z^3+\cdots
$$
$$
=\left(1+{z\over w}+{z^2\over w^2}+\cdots\right)b_0+\left(1+{z\over w}+{z^2\over w^2}+\cdots\right)b_1z+\left(1+{z\over w}+{z^2\over w^2}+\cdots\right)b_2z^2+\cdots
$$
$$
\hskip20pt=\left(1+{z\over w}+{z^2\over w^2}+\cdots\right)\sigma_+(z)b_0=i_{w,z}{w\over w-z}\sigma_+(z)\sigma_-(w)b_0.\hskip141pt \qed
$$
\bclm{\bf Lemma.} {\em For all $n\geq 0$:
\be
\sigma_-(w)\sigma_+(z)b_n=\sigma_+(z)\sigma_-(w)b_n+{1\over w^n}i_{w,z}{z\over w-z}\sigma_+(z)\sigma_-(w)b_0.\label{eq5:20}
\ee
}
\eclm
\proof
First we use formula (\ref{eq5:19}):
\begin{center}
\begin{tabular}{rllr}
$\sigma_-(w)\sigma_+(z)b_n$&$=$&$\sigma_-(w)b_n +\sum_{i\geq 1}\sigma_-(w)b_{n+i}z^i$&(definition of $\sigma_+(z)$)\cr\cr
&$=$&$\displaystyle{\sum}_{i\geq 0}\sigma_i\sigma_-(w)b_nz^i+\displaystyle{\sum}_{i\geq 1}\displaystyle{1\over w^{n+1}}\sigma_-(w)b_{i-1}z^i$&(by eq. (\ref{eq5:19}))\cr\cr
&$=$&$\sigma_+(z)\sigma_-(w)b_n+\displaystyle{z\over w^{n+1}}\sigma_-(w)\sigma_+(z)b_0$&(again definition of \cr &&&$\sigma_+(z)$)\cr\cr
&$=$&$\sigma_+(z)\sigma_-(w)b_n+\displaystyle{z\over w^{n+1}}i_{w,z}{w\over w-z}\sigma_+(z)\sigma_-(w)b_0$&(Lemma \ref{lem:lem56})\cr\cr
&$=$&$\sigma_+(z)\sigma_-(w)b_n+\displaystyle{1\over w^{n}}i_{w,z}\displaystyle{z\over w-z}\sigma_+(z)\sigma_-(w)b_0$.&(simplification)
\end{tabular}
\end{center}
\qed

In particular, for $n=0$:
\begin{eqnarray}
\sigma_-(w)\sigma_+(z)b_0&=&\sigma_+(z)\sigma_-(w)b_0+i_{w,z}{z\over w-z}\sigma_+(z)\sigma_-(w)b_0\cr\cr
&=&i_{w,z}{w\over w-z}\sigma_+(z)\sigma_-(w)b_0.\label{eq5:comb0}
\end{eqnarray}
\bclm{\bf Proposition.} \label{prop5:55}{\em Let $\blamb\in\Pcal_r$. Then:
\[
\sigma_-(w)\sigma_+(z)\wb^{r+1}_\blamb=i_{w,z}{w\over w-z}\sigma_+(z)\sigma_-(w)\wb^{r+1}_\blamb.
\]
}
\eclm
\proof
By using (\ref{eq5:20}) and (\ref{eq5:comb0}):
\begin{eqnarray*}
&&\sigma_-(w)\sigma_+(z)\wb^{r+1}_\blamb\cr\cr
&\hskip-5pt =&\sigma_-(w)\sigma_+(z)\left(b_{r+\lambda_1}\w\cdots\w b_{\lambda_r}\w b_0\right)\cr\cr
&\hskip-5pt=&\sigma_-(w)\sigma_+(z)b_{r+\lambda_1}\w\cdots\w \sigma_-(w)\sigma_+(z)b_{1+\lambda_r}\w \sigma_-(w)\sigma_+(z)b_0\cr\cr
&\hskip-7pt =&\hskip-8pt \bigwedge_{i=1}^{r}\left(\sigma_+(z)\sigma_-(w)b_{r-i+1+\lambda_i}+{z\over w^{r+1+\lambda_1}}i_{w,z}{w\over w-z}\sigma_+(z)\sigma_-(w)b_0\right)\w i_{w,z}{w\over w-z}\sigma_+(z)\sigma_-(w)b_0\cr\cr\cr &\hskip-5pt =&i_{w,z}{w\over w-z}\sigma_+(z)\sigma_-(w)\wb^{r+1}_\blamb.\hskip276pt\qed
\end{eqnarray*}

\bclm{\bf Lemma.} \label{lem5:lemw} For all $\blamb\in\Pcal_r$:
$$
\ovsig_-(w)(\beta_0\lrcorner \wb^{r+1}_\blamb)=w^{-r}\ovsig_+(w)\wb^r_\blamb.
$$
\eclm
\proof If $r=1$, then
$$
\ovsig_-(w)(\beta_0\lrcorner b_{1+\lambda}\w b_0)=w^{-1}(b_\lambda-b_{\lambda+1}w)=w^{-1}\ovsig_+(w)b_\lambda
$$
and the property holds for $r=1$. For $r\geq 1$ it follows by using the fact that $\ovsig_-(w)$ is a HS derivation. In fact
$$
\ovsig_-(w)(\beta_0\lrcorner \wb^{r+1}_\blamb)=\ovsig_-(w)\wb^r_\blamb=[\ovsig_-(w)\bfb]^r_\blamb=w^{-r}\ovsig_+(z)\wb^r_\blamb
$$
\qed

\bclm{\bf Theorem.}\label{thm5:2ndthm} {\em Let
\be
\Gamma_r(z,w):={E_r(w)\over E_r(z)}\sigma_-(w)\ovsig_-(z).
\ee
Then
\be
\Ecal(z,w)=i_{w,z}{w\over w-z}\left({z^r\over w^r}\Gamma_r(z,w)-1\right).
\ee
}
\eclm
\proof
By properly expanding the expression $w{\bm\beta}(w)\lrcorner (\bfb(z)\w \wb^{r}_\blamb)$ and using \ref{sec3:35} one obtains:
\begin{center}
\begin{tabular}{rcll}
$w{\bm\beta}(w)\lrcorner (\bfb(z)\w \wb^{r}_\blamb)$&$=$&$z^r\sigma_-^T(w)
\beta_0\lrcorner \left(\sigma_+(z)\ovsig_-(z)
\wb^{r+1}_\blamb\right)$& (expression of $\bfb(z)$ and ${\bm\beta}(w)$\cr
&&& through Schubert derivations)\cr\cr
&$=$&$z^r\ovsig_-(w)\left(
\beta_0\lrcorner\sigma_-(w)\sigma_+(z)\ovsig_-(z)\wb^{r+1}_\blamb\right)$&(integration by parts (\ref{eq2:iptra}))\cr\cr
&$=$&$\displaystyle{z^r\over w^r}\ovsig_+(w)\sigma_-(w)\sigma_+(z)\ovsig_-(z)\wb^r_\blamb$& (by Lemma \ref{lem5:lemw})\cr\cr
&$=$&$\displaystyle{z^r\over w^r}i_{w,z}{w\over w-z}\ovsig_+(w)\sigma_+(z)\left(\sigma_-(w)\ovsig_-(z)\wb^r_\blamb\right)$& (Proposition \ref{prop5:55})\cr\cr
&$=$&${z^r\over w^r}i_{w,z}\displaystyle{w\over w-z}{E_r(w)\over E_r(z)}\sigma_-(w)\ovsig_-(z)\wb^r_\blamb.$&(invoking $B_r$-module\cr
&&& \,\,\, structure of $\bw^rM$)
\end{tabular}
\end{center}

Therefore
$$
\Ecal(z,w)=i_{w,z}{w\over w-z}-{z^r\over w^r}i_{w,z}{w\over w-z}{E_r(w)\over E_r(z)}\sigma_-(w)\ovsig_-(z)=i_{w,z}{w\over w-z}\cdot\left({z^r\over w^r}\,\Gamma_r(z,w)-1\right),
$$
as desired.\qed
\claim{} In the ring $B_r\otimes_\ZZ\QQ$ one can define the sequence $(x_1,x_2,\ldots)$ related to $(e_1,\ldots, e_r)$ by the relation
$$
\exp\left(\sum_{i\geq 1}x_iz^i\right)E_r(z)=1
$$
In \cite[Theorem 7.1]{pluckercone} it is shown that,  for $r=\infty$:
$$
\ovsig_-(z)=\exp\left(-\sum_{i\geq 1}{1\over iz^i}{\partial\over \partial x_i}\right)\qquad\mathrm{and}\qquad \sigma_-(w)=\exp\left(-\sum_{i\geq 1}{1\over iw^i}{\partial\over \partial x_i}\right).
$$
Thus:
$$
\Gamma_\infty(z,w)={E_\infty(w)\over E_\infty(z)}\sigma_-(w)\ovsig_-(z)=\exp\left(\sum_{i\geq 1}x_i(z^i-w^i)\right)\exp\left(-\sum_{i\geq 1}{1\over i}\left({1\over z^i}-{1\over w^i}\right){\partial\over \partial x_i}\right).
$$
the classical expression of the vertex operator involved in the DJKM bosonic representation of $gl_\infty(\QQ)$.
\section{The $gl_n(\ZZ)$ structure of $B_{r,n}$.}

\claim{} Recall that $M_n:=\bigoplus_{0\leq j<n} \ZZ\cdot b_j$. The submodule of  $\sigma_nM:=\bigoplus_{j\geq n}\ZZ\cdot b_j\subseteq M$ sits into the split exact sequence:
$$
0\sra \sigma_nM\sra M\sra M_n\sra 0,
$$
so that $M_n$ can be identified with the quotient $M/\sigma_nM$. Similarly, the module $\bw^rM_n$ sits into the bottom exact sequence of the following commutative diagram 
$$
\matrix{0&\lra&(h_{n-r+1},\ldots, h_n)&\hookrightarrow&B_r&\stackrel{\pi_{r,n}}\lra&B_{r,n}&\lra&0\cr
&&\Big\downarrow&&\Big\downarrow&&\Big\downarrow&\cr
0&\lra &\bw^{r-1}M\w \sigma_nM&\hookrightarrow&\bw^rM&\stackrel{\pi_{r,n}}\lra& \bw^rM_n&\lra &0} 
$$
whose vertical arrows are multiplication by $\wb^r_0$ and where abusing notation we have denoted by $\pi_{r,n}$ both the canonical projection $B_r\sra B_{r,n}$ ad $\bw^rM\sra \bw^rM_n$.

\bclm{} Let $I_{r,n}$ the ideal $(h_{n-r+1},\ldots,h_n)$. Under the $B_r$-module structure (\ref{eq:Bmodst}) of $\bw^rM$
\be
I_{r,n}\wb^r_0=\sigma_nM\w \bw^{r-1}M.
\ee
\eclm
\proof
Indeed, $I_{r,n}\wb^r_0\subseteq \sigma_nM\w \bw^{r-1}M$, because
$$
h_{n-r+1+j}\wb^r_0=b_{n+j}\w b_{r-2}\w \cdots\w b_0\in \sigma_nM\w \bw^{r-1}M,
$$ 
for 
all $j\geq 0$. Conversely, if $\wb^r_\blamb\in\sigma_nM\w \bw^{r-1}M$, then $
\lambda_1\geq n-r+1$. If $\Delta_\blamb(H_r)\wb_0^r=\wb^r_\blamb$, then $\Delta_\blamb(H_r)$ belongs to the ideal generated by the its first column $(h_{\lambda_1},h_{\lambda_1+1}.\ldots,h_{\lambda_1+r_1})$. Since $\lambda_1\geq n-r+1$, it follows that $(h_{\lambda_1},h_{\lambda_1+1}.\ldots,h_{\lambda_1+r_1})\subseteq I_{r,n}$, i.e. $\wb^r_\blamb\in I_{r,n}\wb^r_0$.\qed

\claim{} Let 
$$
\delta(z,w)_n:=\sum_{0\leq i,j< n}\delta(b_i\otimes \beta_j)z^iw^{-j},
$$
 and define
$$
\Ecal(z,w)_n=\sum_{0\leq i,j<n}b_i\otimes \beta_j\cdot z^iw^{-j}.
$$
as  a map $B_{r,n}\sra B_{r,n}[z,\wi]$ through the equality:
\be
(\Ecal(z,w)_n\Delta_\blamb(H_{r,n}))\left(\wb^r_0+ \sigma_nM\w \bw^{r-1}M\right)=\delta(z,w)_n\left(\wb^r_\blamb+\sigma_nM\w \bw^{r-1}M\right).\ee

\bclm{\bf Theorem.}\label{thm6:mnthm} {\em The $gl_n(\ZZ)$-module structure of $B_{r,n}$ is described by:
\be
\Ecal(z,w)_n\Delta_\blamb(H_{r,n})={z^{r-1}\over w^{r-1}}\pi_{r,n}\left({1\over E_r(z)}\right)\pi_{r,n}\ovsig_-(z)\Delta_\blamb(w^{-\blamb},H_r),\label{eq6:mnthm}
\ee
or, more explicitly:
\be
\Ecal(z,w)_n\Delta_\blamb(H_{r,n})={z^{r-1}\over w^{r-1}}(1+h_1z+\ldots+h_{n-r}z^{n-r})\pi_{r,n}\ovsig_-(z)\Delta_\blamb(w^{-\blamb},H_r).\label{eq6:mnthm1}
\ee
}
\eclm
\proof Since
\[
\delta_n(z,w)\left(\sigma_nM\w \bw^{r-1}M\right)\subseteq  \sigma_nM\w \bw^{r-1}M,
\]
as a simple exercise shows, it follows that
$\Ecal(z,w)_nI_{r,n}\subseteq I_{r,n}$. Therefore:
\begin{eqnarray}
\Ecal(z,w)_n \Delta_\blamb(H_{r,n})\left(\wb^r_0+\bw^{r-1}M\w \sigma_nM\right)&=&\delta_n(z,w)\left(\wb^r_\blamb+\sigma_nM\w \bw^{r-1}M\right)\cr
&=&\delta(z,w)\wb^r_\blamb+\bw^{r-1}M\w \sigma_nM\cr\cr
&=&(\Ecal(z,w) \Delta_\blamb(H_r))\wb^r_0+\sigma_nM \w\bw^{r-1}M,
\end{eqnarray}
i.e., in other words, 
$
\Ecal(z,w)_n\Delta_\blamb(H_{r,n})=\pi_{r,n}\left(\Ecal(z,w) \Delta_\blamb(H_r)\right).
$
Using Theorem \ref{thm:thm43} and the fact that $\pi_{r,n}$ is a epimorphism, one 
finally obtains (\ref{eq6:mnthm}). Equality~(\ref{eq6:mnthm1}) follows by noticing that
$$
\hskip82pt \pi_{r,n}\left({1\over E_r(z)}\right)=\pi_{r,n}(\sum_{i\geq 0}h_iz^i)=1+h_1z+\cdots+h_{n-r}z^{n-r}.\hskip 82pt \qed
$$

\claim{\bf Remark.} It is important to notice that $\pi_{r,n}\circ\ovsig_-(z)\neq 
\ovsig_-(z)\circ \pi_{r,n}$ and that 
$$
\pi_{r,n}\ovsig_-(z)\Delta_\blamb(w^{-\blamb}, H_r)= \ovsig_-(z)\Delta_\blamb(w^{-
\blamb}, \pi_{r,n}H_r)=\ovsig_-(z)\Delta_\blamb(w^{-\blamb},H_{r,n})
$$
only if $n> r-1+\lambda_1$. Thus formula (\ref{eq6:mnthm}) is already in its best 
possible shape.

\bclm{\bf Example.}
  Let us evaluate $\Ecal(z,w)_4\Delta_{(2,2)}(H_{2,4})\in B_{2,4}[z,\wi]$, where $B_{2,
  4}$ is th cohomology (or Chow) ring of the Grassmannian $G(2,4)$. Recall that in this 
  case $h_i=c_i(\Qcal_2)$, the $i$-th Chern class of the universal quotient bundle over 
  it.
According to the recipe, we first compute
$$
\sigma_-(z)\Delta_{(2,2)}(w^{-(2,2)},H_2)=\left|\matrix{w^{-2}&h_1\cr\cr h_3-
\displaystyle{h_2\over z}& h_2-\displaystyle{h_1\over z}}\right|.
$$
Projecting ont $B_{2,4}$ via $\pi_{2,4}$ amounts to set $h_3$ to $0$. Then we muliply 
by $\pi_{2,4}(1/E_2(z))=1+h_1z+h_2z^2$ and by $z/w$.
Finally we obtain:
\begin{eqnarray}
\Ecal(z,w)_4\Delta_{(2,2)}(H_{2,4})&=&{z\over w}(1+h_1z+h_2z^2)\left|\matrix{w^{-2}
&w^{-1}\cr\cr -\displaystyle{h_2\over z}& h_2-\displaystyle{h_1\over z}}\right|=\cr\cr
\cr
&=&h_2{1\over w^2}+ {h_1h_2}{z\over w^2}+h_2^2{z^2\over w^{2}}-h_1{1\over w^3}-(h_1^2-h_2){z\over w^3}+h_2^2{z^2\over w^3}
\end{eqnarray}
so, for instance,
$$
\pmatrix{0&0&0&0\cr 0&1&0&0\cr 0&0&0&0\cr 0&0&0&0}\Delta_{(2,2)}(H_{2,4})=\Ecal_{1,2}\Delta_{(2,2)}(H_{2,4})=h_1h_2
$$
\eclm
\claim{\bf Example.} Recall Example \ref{ex:ex44}, where we have computed
$$
\Ecal_{4,2}e_2=h_1h_3-h_4.
$$
This is zero in $B_{2,4}:=B_2/(h_3,h_4)$. Indeed
\begin{eqnarray*}
\Ecal(z,w)_2e_2&=&\left[{1\over w}(h_1-h_2z)+{1\over w^2}(h_1z-1)\right](1+h_1z+h_2z^2)\cr\cr\cr &=&[e_1+e_2z+(2e_1^2e_2-e_1^4-e_2^2z^3]w^{-1} +[-1+e_2z^2+(e_1^3-e_1e_2)z^3]w^{-2}.
\end{eqnarray*}

\bibliographystyle{amsplain}
\providecommand{\bysame}{\leavevmode\hbox to3em{\hrulefill}\thinspace}
\providecommand{\MR}{\relax\ifhmode\unskip\space\fi MR }
\providecommand{\MRhref}[2]{%
  \href{http://www.ams.org/mathscinet-getitem?mr=#1}{#2}
}
\providecommand{\href}[2]{#2}

\parbox[t]{3in}{{\rm Letterio~Gatto}\\
{\tt \href{mailto:letterio.gatto@polito.it}{letterio.gatto@polito.it}}\\
{\it Dipartimento~di~Scienze~Matematiche}\\
{\it Politecnico di Torino}\\
{\it ITALY}} \hspace{1.5cm}
\parbox[t]{2.5in}{{\rm Parham~Salehyan}\\
{ \href{mailto:p.salehyan@unesp.br}{\tt{p.salehyan@unesp.br}}}\\
{\it Ibilce UNESP}\\
{\it Campus de S\~ao Jos\'e do Rio Preto, SP}\\
{\it BRAZIL}}

\end{document}